\documentclass[12pt,a4paper,reqno]{amsart}
\usepackage[ansinew]{inputenc}
\usepackage{amsthm,amsmath,amssymb}
\usepackage{enumerate,verbatim}
\RequirePackage[colorlinks,citecolor=blue,urlcolor=blue]{hyperref}

\newcommand\NoBlackBoxes{\global\overfullrule0pt}
\NoBlackBoxes
\theoremstyle{plain} 
\def\4{\kern1pt}

\def\6{\vphantom0}

\def\8{\kern-10pt}
\def\7#1{_{(#1)}}

\makeatletter
\let\serieslogo@\relax
\let\@setcopyright\relax

\def\speciallabelmark#1{\def\@currentlabel{#1}}
\makeatother

\begin{document}

\def\ffrac#1#2{\raise.5pt\hbox{\small$\4\displaystyle\frac{\,#1\,}{\,#2\,}\4$}}
\def\ovln#1{\,{\overline{\!#1}}}
\def\ve{\varepsilon}
\def\kar{\beta_r}
\def\theequation{\thesection.\arabic{equation}}
\def\E{{\mathbb E}}
\def\R{{\mathbb R}}
\def\C{{\mathbb C}}
\def\Z{{\mathbb Z}}
\def\P{{\mathbb P}}
\def\H{{\rm H}}
\def\Im{{\rm Im}}
\def\Tr{{\rm Tr}}

\def\k{{\kappa}}
\def\M{{\cal M}}
\def\Var{{\rm Var}}
\def\Ent{{\rm Ent}}
\def\O{{\rm Osc}_\mu}

\def\ep{\varepsilon}
\def\phi{\varphi}
\def\F{{\cal F}}
\def\L{{\cal L}}

\def\s{{\mathfrak s}}

\def\be{\begin{equation}}
\def\en{\end{equation}}
\def\bee{\begin{eqnarray*}}
\def\ene{\end{eqnarray*}}


\title[CLT for R\'enyi Divergence]{CENTRAL LIMIT THEOREM FOR R\'ENYI \\
DIVERGENCE OF INFINITE ORDER}


\author{S. G. Bobkov$^{1}$}
\thanks{1) 
School of Mathematics, University of Minnesota, Minneapolis, MN, USA,
bobkov@math.umn.edu.
}

\author{F. G\"otze$^{2}$}
\thanks{2) Faculty of Mathematics,
Bielefeld University, Germany,
goetze@math-uni.bielefeld.de.
}






\begin{abstract}
For normalized sums $Z_n$ of i.i.d. random variables, we explore necessary
and sufficient conditions which guarantee the normal approximation with 
respect to the R\'enyi divergence of infinite order. In terms of densities 
$p_n$ of $Z_n$, this is a strengthened variant of the local limit theorem taking the form 
$\sup_x (p_n(x) - \varphi(x))/\varphi(x) \rightarrow 0$ as $n \rightarrow \infty$.
\end{abstract}




\maketitle

\section{{\bf Introduction. Strict Subgaussianity}}
\setcounter{equation}{0}

\vskip2mm
\noindent
Let $X$ be a random variable with density $p$. The R\'enyi divergence 
of order $\alpha > 0$, or the relative $\alpha$-entropy of its distribution 
with respect to the standard normal law with density 
$\varphi(x) = \frac{1}{\sqrt{2\pi}}\,\exp(-x^2/2)$ is given by
\be
D_\alpha(p||\varphi) = \frac{1}{\alpha - 1} \log 
\int_{-\infty}^\infty \Big(\frac{p}{\varphi}\Big)^\alpha\,\varphi\, dx.
\en
A closely related functional is the Tsallis distance
\be
T_\alpha(p||\varphi) = \frac{1}{\alpha - 1} \bigg[\int_{-\infty}^\infty 
\Big(\frac{p}{\varphi}\Big)^\alpha\,\varphi\, dx - 1\bigg].
\en
Since
$T_\alpha = \frac{1}{\alpha - 1}\,[e^{(\alpha-1)\,D_\alpha} - 1]$,
both distances are of a similar order, when they are small.
Hence, approximation problems in $D_\alpha$ and $T_\alpha$
are equivalent. Moreover, as the function 
$\alpha \rightarrow D_\alpha$ is non-decreasing, the convergence in 
$D_\alpha$ is getting stronger for growing indexes $\alpha$.

Let us recall that, for the region $0 < \alpha < 1$, $D_\alpha$ 
is topologically equivalent to the total variation distance 
between the distribution of $X$ and the standard normal law.
For $\alpha=1$, we obtain the Kullback-Leibler distance
$$
D(p||\varphi) = \int_{-\infty}^\infty p\log \frac{p}{\varphi}\, dx,
$$
also called the informational divergence or the relative entropy.
It is finite, if and only if $X$ has a finite second moment and
finite Shannon's entropy. But,
the range $\alpha>1$ leads to much stronger R\'enyi/Tsallis 
distances. For example, the finiteness of $D_\alpha(p||\varphi)$
requires that $X$ is subgaussian, i.e. the moments $\E\,e^{cX^2}$
should be finite for small $c>0$. One important 
particular case $\alpha = 2$ in this hierarchy corresponds
to the Pearson $\chi^2$-distance $T_2=\chi^2$. 
For various properties and applications of these distances, we refer 
an interested reader to \cite{LC}, \cite{S}, \cite{D-C-T}, 
\cite{J}, \cite{VE-H}, \cite{B-C-G3}.

The study of the convergence in the central limit theorem (CLT) with 
respect to $D_\alpha$ and the associated problem of Berry-Esseen 
bounds have a long and rich history. Let us remind several results 
in this direction about the classical model of normalized sums
$$
Z_n = (X_1 + \dots + X_n)/\sqrt{n}
$$
of i.i.d. random variables $(X_k)_{k \geq 1}$.  
We will treat them as independent copies of a random variable $X$,
assuming that it has mean zero and variance one.

The convergence $D_\alpha(p_n||\varphi)\rightarrow 0$ as $n \rightarrow \infty$ 
holds true for $0 < \alpha < 1$, as long as $Z_n$ have densities 
$p_n$ for large $n$. This is due to the corresponding result by Prokhorov \cite{Pr}
about the total variation distance. The stronger property $D(p_n||\varphi)\rightarrow 0$ 
in terms of relative entropy was studied by Barron \cite{Bar} who showed that 
the condition $D(p_n||\varphi) < \infty$ for some $n$ is necessary and sufficient for 
the entropic CLT. The asymptotic behavior of such distances under 
higher order moment assumptions, including Edgeworth-type expansions in powers 
of $1/n$, has been studied in \cite{B-C-G1}. It is worthwhile mentioning that this 
convergence is monotone with respect to $n$, cf.\,Artstein, Ball, Barthe and 
Naor \cite{A-B-B-N1} and Madiman and Barron \cite{M-B}. See also \cite{A-B-B-N2}
and \cite{B-C-G2} for various entropic bounds in the non-i.i.d.~case.

The range $\alpha>1$ was treated 
in detail in \cite{B-C-G3}. It was shown there that $D_\alpha(p_n||\varphi)\to 0$ 
as $n\to\infty$, if and only if $D_\alpha(p_n||\varphi)$ is finite for some $n$, 
and if $X$ admits the following subgaussian bound on the Laplace transform
\be
\E\,e^{tX} < e^{\alpha^* t^2/2}, \quad t \in \R \ (t\ne 0),
\en
where $\alpha^* = \frac{\alpha}{\alpha - 1}$. 
In that case, we have an equivalence
$D_\alpha \sim T_\alpha \sim \frac{\alpha}{2} \chi^2$.
These results have been extended to the multidimensional setting as well.

For indexes $\alpha \rightarrow \infty$ in (1.3), we
arrive at the following characterization:

\vskip5mm
{\bf Theorem 1.1.} {\sl Assume that $D_\alpha(p_n||\varphi) <\infty$
for every $\alpha>1$ with some $n = n_\alpha$. For
the convergence $D_\alpha(p_n||\varphi) \rightarrow 0$ for all $\alpha$,  
it is necessary and sufficient that $\E \exp\{tX\} \leq \exp\{t^2/2\}$
for all $t \in \R$.
}

\vskip5mm
The last inequality describes an interesting class of probability distributions
which appear naturally in many mathematical problems. More generally, one says 
that a random variable $X$ with mean zero is strictly subgaussian, or its
distribution is strictly subgaussian (regardless of whether or
not it has density), if the inequality
\be
\E\,e^{tX} \leq e^{\sigma^2 t^2/2}, \quad t \in \R,
\en
holds with constant $\sigma^2 = \Var(X)$ which is then best possible. 
Note that, when saying that $X$ is subgaussian, one means that (1.4) holds 
with some~$\sigma^2$.

This class was apparently first introduced in an explicit form by
Buldygin and Kozachenko in \cite{Bu-K1} under the name ``strongly subgaussian" 
and then analyzed in more details in their book \cite{Bu-K2}.
Recent investigations include the work by Arbel, Marchal and Nguyen \cite{A-M-N} 
providing some examples and properties and by Guionnet and Husson \cite{G-H}.
In the latter paper, (1.4) appears as a condition for the validity of large deviation 
principles for the largest eigenvalue 
of Wigner matrices with the same rate function as in the case of Gaussian entries. 

A simple sufficient condition for the strict subgaussianity was given by 
Newman in terms of location of zeros of the characteristic function
$f(z) = \E\,e^{izX}$, $z \in \C$ (which is extended, by the subgaussian property, 
from the real line to the complex plane as an entire function of order at most 2). 
As was stated in \cite{N1}, $X$ is strictly subgaussian, as long as $f(z)$ has 
only real zeros in $\C$ (a detailed proof was later given in \cite{Bu-K2}).
Such probability distributions form an important class denoted by $\mathfrak L$, 
introduced 
and studied by Newman in the mid 1970's in connection with the Lee-Yang property 
which naturally arises in the context of ferromagnetic Ising models, 
cf. \cite{N1,N2,N3,N-W}. This class is rather rich; it is closed under infinite
convergent convolutions and under weak limits. For example, it includes 
Bernoulli convolutions and hence convolutions of uniform distributions
on bounded symmetric intervals.

Some classes of strictly subgaussian distributions outside $\mathfrak L$
have been recently discussed in \cite{B-C-G5}. It was shown that (1.4) continues 
to hold under the weaker requirement that all zeros of $f(z)$ 
with ${\rm Re}(z) > 0$ lie in the cone
$|{\rm Arg}(z)| \leq \frac{\pi}{8}$ (which is sharp when $f$ has only
one zero in the positive octant). In that case, if $X$ 
is not normal, the inequality (1.4) may be sharpened 
as follows: For any $t_0 > 0$, there is $c = c(t_0)$,
$0 < c < \sigma^2 = \Var(X)$, such~that
\be
\E\,e^{tX} \leq e^{c t^2/2}, \quad |t| \geq t_0.
\en

In general, this separation-type property is however not necessary
for the strict subgaussianity. It turns out that there exists a large class of
strictly subgaussian distributions with mean zero and variance one, 
for which the Laplace transform has the form
$$
\E\,e^{tX} = \Psi(t)\, e^{-t^2/2}, \quad t \in \R,
$$
where $\Psi(t)$ is a {\it periodic} function with some period $h>0$ and such that
$\Psi(t) \leq 1$ for all $t \in \R$. Hence $\Psi(kh) = 1$ for all $k \in \Z$,
so that (1.4) becomes an equality for infinitely many points $t$.

\section{{\bf Main Results for the Convergence in $D_\infty$}}
\setcounter{equation}{0}

Thus, the strict subgaussianity appears as 
a necessary condition for the convergence in all $D_\alpha$ and 
therefore in $D_\infty$, which according to (1.1) is given by the limit
$$
D_\infty(p||\varphi) = \lim_{\alpha \rightarrow \infty} D_\alpha(p||\varphi) =
{\rm ess\,\sup}_x\, \log(p(x)/\varphi(x)).
$$
Although the Tsallis distance of infinite order may not be
defined similarly as a limit of (1.2), we make the convention that
$$
T_\infty(p||\varphi) = {\rm ess\,\sup}_x\, \frac{p(x) - \varphi(x)}{\varphi(x)}.
$$
Then $T_\infty = e^{D_\infty} - 1$
like for the Tsallis distance of finite order, so that convergence in 
$D_\infty$ and $T_\infty$ are equivalent.
In particular, in the setting of the normalized sums $Z_n$,
the CLT $D_\infty(p_n||\varphi) \rightarrow 0$ is equivalent
to the assertion that $Z_n$ have densities $p_n$ such that
\be
\sup_x \, \frac{p_n(x) - \varphi(x)}{\varphi(x)} \rightarrow 0 \quad
{\rm as} \ n \rightarrow \infty.
\en

The purpose of this paper is to give necessary and sufficient conditions
for this variant of the CLT in terms of the Laplace transform
$L(t) = \E\,e^{tX}$. Consider the log-Laplace transform
$K(t) = \log L(t)$ (which is a convex function) and the associated function
$$
A(t) = \frac{1}{2} t^2 - K(t), \quad t \in \R.
$$
As before,
suppose that $(X_k)_{k \geq 1}$ are independent copies of the random
variable $X$ with $\E X = 0$ and $\Var(X) = 1$. We assume that:

\vskip2mm
1) $Z_n$ has density $p_n$ with $T_\infty(p_n||\varphi) < \infty$ for some $n=n_0$;

2) $X$ is strictly subgaussian, that is, $A(t) \geq 0$ for all $t \in \R$.

\vskip5mm
{\bf Theorem 2.1.} {\sl For the convergence $T_\infty(p_n||\varphi) \rightarrow 0$,
it is necessary and sufficient that the following two conditions are fulfilled:

\vskip2mm
$a)$ \ $A''(t) = 0$ for every point $t \in \R$ such that $A(t) = 0$;

$b)$ \ $\limsup_{k \rightarrow \infty} A''(t_k) \leq 0$ for every sequence
$t_k \rightarrow \pm \infty$ such that $A(t_k) \rightarrow 0$ as $k \rightarrow \infty$.
}

\vskip2mm
The conditions $a)-b)$ may be combined as
$\lim_{A(t) \rightarrow 0} \max(A''(t),0) = 0$,
which is kind of continuity of $A''$ with respect to $A$.

Under the separation property (1.5), the condition $b)$ is fulfilled automatically,
while the equation $A(t) = 0$ has only one solution $t=0$. But near zero,
due to the strict subgaussianity, $A(t) = O(t^4)$ and $A''(t) = O(t^2)$.
Hence, the condition $a)$ is fulfilled as well, and we obtain the CLT with
respect to $D_\infty$. In particular, it is applicable to the class $\mathfrak L$ 
of Newman described above. In fact, for this conclusion, (1.5) may further be weakened to
\be
\sup_{|t| \geq t_0} \big[e^{-t^2/2}\, \E\,e^{tX} \big] < 1 \quad {\rm for \ all} \ t_0>0.
\en
In this case one can additionally explore the rate of convergence.

\vskip5mm
{\bf Theorem 2.2.} {\sl Let $X$ be a non-normal random variable
with $\Var(X) = 1$ satisfying $(2.2)$. If $T_\infty(p_n||\varphi) < \infty$ 
for some $n$, then
\be
T_\infty(p_n||\varphi) = O\Big(\frac{1}{n}\,(\log n)^3\Big)
\quad {\sl as} \ n \rightarrow \infty.
\en
}

Furthermore, specializing Theorem 2.1 to the case where the Laplace transform
contains a periodic component, we have:

\vskip5mm
{\bf Theorem 2.3.} {\sl Suppose that the function $\Psi(t) = L(t)\,e^{-t^2/2}$
is $h$-periodic for a smallest value $h>0$. For the convergence
$T_\infty(p_n||\varphi) \rightarrow 0$ as $n \rightarrow \infty$,
it is necessary and sufficient that, for every $0 < t < h$,
\be
\Psi(t) = 1 \, \Rightarrow \, \Psi''(t) = 0.
\en
Moreover, if the equation $\Psi(t) = 1$ has no solution in $0 < t < h$, then
the relation $(2.3)$ about the rate of convergence continues to hold.
}

\vskip5mm
For an illustration  (cf. Section 9 for more details) , consider random variables $X$ with
$\Psi(t) = 1 - c\sin^4 t$, where the parameter $c>0$ is small enough.
In this case, $\Psi(t)$ is $\pi$-periodic and all conditions in Theorem 2.1
are fulfilled. Hence the CLT for $T_\infty$ does hold with rate as in (2.3).
On the other hand, in a similar $\pi$-periodic example
$$
\Psi(t) = 1 - c\,(1 - 4\sin^2 t)^2 \, \sin^4 t,
$$
the condition (2.4) is violated at the point $t = \pi/6$, so there is no CLT.
Thus, the continuity condition of $A''$ with respect $A$ in Theorem 2.1
may or may not be fulfilled in general in the class of strictly subgaussian distributions.

Returning to the convergence property (2.1), it should be emphasized 
that it is not possible to put the absolute value sign in the numerator
(this will be clarified in Section~4). The situation is of course different, 
when one considers the supremum over bounded increasing intervals. 
For example, under suitable moment assumptions (cf. \cite{Pe1}, \cite{Pe2}), 
it follows from Edgeworth expansions for densities that
$$
\sup_{|x| \leq c\sqrt{\log n}} \, \frac{|p_n(x) - \varphi(x)|}{\varphi(x)}
\rightarrow 0 \quad {\rm as} \ n \rightarrow \infty.
$$

The proof of Theorem 2.1 is given in Section 8, with preliminary developments 
in Sections 3-7. Its application to the periodic case is discussed in Section 9.
What is unusual in our approach is that the proof does not use in essence
the tools from Complex Analysis (as one ingredient, we establish a uniform
local limit theorem for bounded densities with a quantitative error term).
However, in the study of rates
of convergence with respect to $T_\infty$, we employ an old result 
by Richter \cite{Ri} in a certain refined form on the asymptotic behavior
of ratios $p_n(x)/\varphi(x)$. This result is discussed in 
Section 10, where we also include the proof of Theorem 2.2 and
Theorem 2.3 (for the rate of convergence). In the last
section, we describe several examples illustrating applicability
of Theorem 2.2.

\section{{\bf Semigroup of Shifted Distributions (Esscher Transform)}}
\setcounter{equation}{0}

Let $X$ be a subgaussian random variable with density $p$. 
Then, the Laplace transform, or the moment generating function
$$
(L p)(t) = L(t) = \E\,e^{tX} = \int_{-\infty}^\infty e^{tx} p(x)\,dx
$$
is finite for all complex numbers $t$ and represents an entire function
in the complex plane. Hence the log-Laplace transform
$$
(K p)(t) = K(t) = \log L(t) = \log \E\,e^{tX}, \quad t \in \R,
$$
represents a convex, $C^\infty$-smooth function on the real line.

\vskip5mm
{\bf Definition 3.1.} Introduce the family of probability densities
\be
Q_h p(x) = \frac{1}{L(h)}\,e^{hx} p(x), \quad x \in \R,
\en
with parameter $h \in \R$. We call the distribution with this density
the shifted distribution of $X$ at step $h$.

\vskip2mm
The early history of this well-known and popular transform goes back 
to 1930's. In actuarial science, following Esscher \cite{E}, the density 
$Q_h p$ is commonly called the Esscher transform of $p$. Other names 
``conjugate distribution laws", ``the family of distribution laws conjugate to a system" 
were used by Khinchin \cite{K} in the framework of statistical mechanics.
See also Daniels \cite{D} who applied this transform to develop
asymptotic expansions for densities. In this paper, we prefer to use
a different terminology as in Definition 3.1 in order to emphasize the following 
important fact: For the standard normal density $\varphi(x)$, 
the shifted normal law has density $Q_h \varphi(x) = \varphi(x+h)$.

A remarkable property of the transform (2.1) is the semi-group property
$$
Q_{h_1}(Q_{h_2} p) = Q_{h_1 + h_2} p, \quad h_1,h_2 \in \R.
$$

Let us also mention how this transform acts under rescaling. Given $\lambda>0$,
the random variable $\lambda X$ has density 
$p_\lambda(x) = \frac{1}{\lambda}\,p(\frac{x}{\lambda})$ with Laplace transform
$(L p_\lambda)(t) = L(\lambda t)$. Hence
$$
Q_h p_\lambda(x) = \frac{1}{(L p_\lambda)(h)}\,e^{hx} p_\lambda(x) = 
\frac{1}{\lambda}\, (Q_{\lambda h} p)\Big(\frac{x}{\lambda}\Big).
$$
This identity implies that the maximum-of-density functional
$$
M(X) = M(p) = {\rm ess\,sup}_x \, p(x)
$$ 
satisfies
\be
M(Q_h p_\lambda) = \frac{1}{\lambda}\, M(Q_{\lambda h} p).
\en

The transform $Q_h$ is also multiplicative with respect to convolutions.

\vskip5mm
{\bf Proposition 3.2.} {\sl If independent subgaussian random
variables have densities $p_1,\dots,p_n$, then
for the convolution $p = p_1 * \dots * p_n$, we have
\be
Q_h p = Q_h p_1 * \dots * Q_h p_n.
\en
}

{\bf Proof.} It is sufficient to compare the Laplace transforms  of both
sides in (3.2). The Laplace transform of $p$ is given by
$
L p(t) = (L p_1)(t) \dots (L p_n)(t).
$
Hence, the Laplace transform of $Q_h p$ is given by
\bee
(L Q_h p)(t) 
 & = & 
\int_{-\infty}^\infty e^{tx}\,Q_h p(x)\,dx \, = \,
\frac{1}{(L p)(t)} \int_{-\infty}^\infty e^{(t + h) x} p(x)\,dx \\
 & = &
\frac{(L p)(t+h)}{(L p)(t)} \, = \, \prod_{k=1}^n \frac{(L p_k)(t+h)}{(L p_k)(t)} \, = \, 
\prod_{k=1}^n (L Q_h p_k)(t).
\ene
\qed

The formula (3.1) in Definition 3.1 may be written equivalently as
$$
p(x) = L(h) e^{-xh}\, Q_h p(x) = e^{-xh + K(h)}\, Q_h p(x),
$$
or
$$
\frac{p(x)}{\varphi(x)} = 
\sqrt{2\pi}\,e^{\frac{1}{2}\, (x-h)^2 - \frac{1}{2} h^2 + K(h)}\, Q_h p(x).
$$
Introduce the function
\be
(A p)(h) =  A(h) = \frac{1}{2} h^2 - K(h),
\en
which allows to reformulate strict subgaussianity  via the inequality $A(h) \geq 0$
for all $h$ (under the assumptions $\E X = 0$ and $\Var(X) = 1$). Thus,
\be
\frac{p(x)}{\varphi(x)} = \sqrt{2\pi}\,e^{\frac{1}{2}\, (x-h)^2 - A(h)}\, Q_h p(x).
\en

We will use this representation to bound the ratio on the left-hand side for 
the densities $p_n$ of the normalized sums
\be
Z_n = (X_1 + \dots + X_n)/\sqrt{n}
\en
of independent copies of the random variable $X$ with density $p$. In order 
to apply (3.5) to $p_n$ instead of $p$, put $x_n = x\sqrt{n}$, $h_n = h\sqrt{n}$.
Note that in terms of $L = L p$, $K = K p$ and $A = A p$, we may write
$$
(L p_n)(t) = L(t/\sqrt{n})^n = e^{n K(t/\sqrt{n})}, \qquad (K p_n)(t) = n K(t/\sqrt{n}),
$$
$$
(A p_n)(h_n) \, = \, \frac{1}{2} h_n^2 - (Kp_n)(h_n) \, = \,
\frac{n}{2} h^2 - n K(h) \, = \, n A(h).
$$
Therefore, the definition (3.5) being applied with $(x_n,h_n)$ becomes:

\vskip5mm
{\bf Proposition 3.3.} {\sl Putting $x_n = x\sqrt{n}$, $h_n = h\sqrt{n}$
$(x,h \in \R)$, we have
\be
\frac{p_n(x\sqrt{n})}{\varphi(x\sqrt{n})} = 
\sqrt{2\pi}\,e^{\frac{n}{2}\, (x - h)^2 - n A(h)}\, Q_{h_n} p_n(x_n).
\en
}

This equality becomes useful, if we are able to bound the factor $Q_{h_n} p_n(x_n)$ 
uniformly over all $x$ for a fixed value of $h$ as stated in the following Corollary.

\vskip5mm
{\bf Corollary 3.4.} {\sl For all $x,h \in \R$,
\be
\frac{p_n(x\sqrt{n})}{\varphi(x\sqrt{n})} \leq 
\sqrt{2\pi}\,e^{\frac{n}{2}\, (x - h)^2 - n A(h)}\, M(Q_{h\sqrt{n}}\, p_n).
\en
}

{\bf Remark 3.5.} Since the function $K$ is convex, it follows from the
definition (3.4) that $A''(h) \leq 1$ for all $h \in \R$. As a consequence,
this function satisfies a differential inequality
\be
A'(h)^2 \leq 2 A(h), \quad h \in \R,
\en
if $A(h) \geq 0$ for all $h \in \R$. For a short proof (proposed by a
referee), one may apply the Taylor formula
\bee
0 \, \leq \, A(h+x) 
 & = &
A(h) + A'(h)x + \frac{1}{2} A''(h_1) x^2 \\
 & \leq &
A(h) + A'(h)x + \frac{1}{2}\,x^2, \quad x \in \R,
\ene
holding for some point $h_1$ in the segment with endpoints $h$ and 
$h+x$. Minimizing the right-hand side over all $x$ leads to (3.9).

\vskip10mm
\section{{\bf Maximum of Shifted Densities}}
\setcounter{equation}{0}

In order to bound the last term in (3.8), suppose that the distribution of $X$ 
has a finite R\'enyi distance of infinite order to the standard normal law. 
This means that the density of $X$ admits a pointwise upper bound
\be
p(x) \leq c \varphi(x), \quad x \in \R \ ({\rm a.e.})
\en
for some constant $c$. Note that its optimal value
is $c = 1 + T_\infty(p||\varphi)$. 
In that case, one may control the maximum 
$$
M(Q_h p) = {\rm ess\,sup}_x\, Q_h p(x)
$$
of densities of shifted distributions. Indeed, (4.1) implies that, for any $x \in \R$,
$$
Q_h p(x) = 
\frac{1}{L(h)}\,e^{xh} p(x) \leq \frac{c\, e^{xh - x^2/2}}{L(h)\sqrt{2\pi}} \leq
\frac{c\, e^{h^2/2}}{L(h)\sqrt{2\pi}} = \frac{c}{\sqrt{2\pi}}\, e^{A(h)},
$$
where $L = Lp$ and $A = Ap$. Thus,
\be
M(Q_h p) \leq \frac{c}{\sqrt{2\pi}}\, e^{A(h)}.
\en

However, it is useless to apply this bound directly to $p_n$ for normalized sums 
$Z_n$ as in (3.6), since then the right-hand side of (4.2) will contain 
the parameter $c_n = 1 + T_\infty(p_n||\varphi)$. Instead, we use
a semi-additive property of the maximum-of-density functional, which indicates that
$$
M(X_1 + \dots + X_n)^{-2} \geq \frac{1}{2} \sum_{k=1}^n M(X_k)^{-2}
$$
for all independent random variables $X_k$ having bounded densities,
cf. \cite{B-C}, p. 105, or \cite{B-C-G4}, p. 142. If all $X_k$ are identically distributed and have density $p$, 
this relation yields
$$
M(p^{*n}) \leq \sqrt{2/n}\,M(p)
$$
for the convolution $n$-th power of $p$. Applying Proposition 3.2 together with (4.2), 
we then have
$$
M(Q_h p^{*n}) \leq 
\sqrt{2/n}\,M(Q_h p) \leq \sqrt{2/n}\,\frac{c}{\sqrt{2\pi}}\, e^{A(h)}.
$$
Now, since $p^{*n}(x) = \frac{1}{\lambda}\,p_n(\frac{x}{\lambda})$ with
$\lambda = \sqrt{n}$, one may apply the identity (3.2):
$$
M(Q_h p^{*n}) = \frac{1}{\sqrt{n}}\, M(Q_{h\sqrt{n}}\, p_n).
$$
Hence
$$
M(Q_{h\sqrt{n}}\, p_n) \leq \frac{c}{\sqrt{\pi}}\, e^{A(h)}.
$$
Now return to Corollary 3.4 and apply this bound to get that
$$
\frac{p_n(x\sqrt{n})}{\varphi(x\sqrt{n})} \leq 
c\sqrt{2}\,e^{\frac{n}{2}\, (x - h)^2 - (n-1) A(h)},
$$
recalling that $c = 1 + T_\infty(p||\varphi)$.
In particular, with $h=x$ this yields:

\vskip5mm
{\bf Proposition 4.1.} {\sl Let $p_n$ denote the density of
$Z_n$ constructed for $n$ independent copies of a subgaussian random variable
$X$ whose density $p$ has finite R\'enyi distance of infinite order to the standard 
normal law. Then, for almost all $x \in \R$,
\be
\frac{p_n(x\sqrt{n})}{\varphi(x\sqrt{n})} \leq 
c\sqrt{2}\,e^{- (n-1) A(x)}.
\en
}

\vskip2mm
{\bf Corollary 4.2.} {\sl If additionally $\E X = 0$, $\Var(X) = 1$, and $X$ is
strictly subgaussian, then
$$
T_\infty(p_n||\varphi) \leq \sqrt{2}\,(1 + T_\infty(p||\varphi)) - 1.
$$
}

\vskip2mm
Thus, the finiteness of the Tsallis distance $T_\infty(p||\varphi)$ for a strictly 
subgaussian random variable $X$ with density $p$ ensures the boundedness
of $T_\infty(p_n||\varphi)$ for all normalized sums $Z_n$.

If $A(x)$ is bounded away from zero, the inequality (4.3) shows that
$p_n(x\sqrt{n})/\varphi(x\sqrt{n})$ is exponentially small for growing $n$.
In particular, this holds for any non-normal random variable $X$ satisfying 
the separation property (2.2). Then we immediately obtain:

\vskip2mm
{\bf Corollary 4.3.} {\sl Suppose that $X$ has a density $p$ with finite 
$T_\infty(p||\varphi)$. Under the condition $(2.2)$, for any $\tau_0 > 0$, 
there exist $A > 0$ and $\delta \in (0,1)$ such that the densities $p_n$ of $Z_n$ satisfy
\be
p_n(x) \leq A\delta^n \varphi(x), \quad
|x| \geq \tau_0\sqrt{n}.
\en
}

In particular,
$$
\liminf_{n \rightarrow \infty} \, \sup_{x \in \R} \,
\frac{|p_n(x) - \varphi(x)|}{\varphi(x)} \geq 1.
$$
Therefore, one can not hope to strengthen the Tsallis distance by
introducing a modulus sign in the definition of the distance.

Since (2.2) does not need be true in general, Proposition 4.1
will be applied outside the set of points where $A(x)$ is bounded away from zero.
More precisely, for a parameter $a>0$ and $n \geq 2$, define the critical zone
\be
A_n(a) = \{h>0: A(h) \leq a/(n-1)\}.
\en
From (4.3), it follows that
\be
\frac{p_n(x\sqrt{n})}{\varphi(x\sqrt{n})} \leq c\sqrt{2}\,e^{-a}, \quad
x \notin A_n(a).
\en
If $a$ is large, this bound may be used in the proof of the CLT with
respect to the distance $T_\infty$ restricted to the complement of the critical zone.
As for this zone, the bound (4.3) is not appropriate, and we need to return
to the basic representation from Proposition 3.2. To study the last term
$Q_{h_n} p_n(x_n)$ in (3.7) for $x \in A_n(a)$, one may apply a variant
of the local limit theorem, using the property that the density $Q_{h_n} p_n$ has
a convolution structure. 
However, in order to justify this application, we should first explore
the behavior of moments of densities participating in the convolution.

\vskip10mm
\section{{\bf Moments of Shifted Distributions}}
\setcounter{equation}{0}

For a subgaussian random variable $X$ with density $p$, denote by 
$X(h)$ a random variable with density $Q_h p$ ($h \in \R$). It is
subgaussian, and its Laplace and log-Laplace transforms are given by
\be
L_h(t) \equiv \E\,e^{t X(h)} = \frac{L(t+h)}{L(h)}, \quad
K_h(t) \equiv \log L_h(t) = K(t+h) - K(h).
\en

Furthermore, it has mean and variance
\bee
m_h 
 & \equiv & \E X(h) \, = \, \frac{L'(h)}{L(h)} \, = \, K'(h), \\
\sigma_h^2 
 & \equiv &
\Var(X(h)) \, = \, \frac{L''(h) - L'(h)^2}{L(h)^2} \, = \, K''(h).
\ene
The last equality shows that necessarily $K''(h) > 0$ for all $h \in \R$.
Indeed, otherwise the random variable $X(h)$ would be a constant a.s.

The question of how to bound the standard deviation
$\sigma_h$ from below relies upon certain fine properties of the density $p$
and the behavior of the function 
$$
A(h) = \frac{1}{2} h^2 - K(h),
$$
introduced in (3.4). As before, suppose that the distribution of $X$ has 
finite R\'enyi distance of infinite order to the standard normal law, so that
\be
p(x) \leq c \varphi(x), \quad x \in \R,
\en
with $c = 1 + T_\infty(p||\varphi)$. Then one may control the maximum 
$M(X(h)) = {\rm ess\,sup}_x\, p_h(x)$
of densities of shifted distributions, using (4.2):
$$
Q_h p(x) \leq  \frac{c}{\sqrt{2\pi}}\, e^{A(h)}.
$$

For a lower bound, we employ a well-known general relation 
$$
M(\xi)^2 \, \Var(\xi) \geq \frac{1}{12}
$$ 
(where the equality is attained for the uniform distribution on 
a bounded interval). Let us provide the following simple argument, assuming 
without loss of generality that a random variable $\xi$ has finite variance
and a density with $M(\xi) = 1$. Then, the function
$H(x) = \P\{|\xi - \E \xi| \geq x\}$ is absolutely continuous, and its
Radon-Nikodym derivative satisfies $H'(x) \geq -2$ a.e. in $x>0$. Since 
$H(0) = 1$, we get $H(x) \geq 1-2x$ for all $x \geq 0$ and therefore
$$
\Var(\xi) \, = \, 
2 \int_0^\infty x H(x)\,dx \, \geq \,
2 \int_0^{1/2} x (1 - 2x)\,dx \, = \, \frac{1}{12}.
$$
Applying this to $\xi = X(h)$ and combining the
two bounds, we obtain that
$$
\frac{1}{\sqrt{12}} \leq M(X(h)) \sigma_h  \leq 
\frac{c\sigma_h}{\sqrt{2\pi}}\, e^{A(h)}.
$$
Thus we arrive at:
\vskip2mm
{\bf Lemma 5.1.} {\sl Under the condition $(5.2)$, for all $h \in \R$,
\be
\sigma_h \geq \sqrt{\frac{\pi}{6c^2}}\, e^{-A(h)}.
\en
}

Since $\sigma_h>0$, one may consider the normalized random variables
\be
\hat X(h) = \frac{X(h) - \E X(h)}{\sqrt{\Var(X(h))}} = \frac{X(h) - m_h}{\sigma_h}.
\en
By (5.1), they have the moment generating function
$$
\E\,e^{t \hat X(h)} \, = \,
\E\,\exp\Big\{\frac{t}{\sigma_h}\,(X(h) - m_h)\Big\} \, = \,
\exp\Big\{-\frac{t}{\sigma_h}\,K'(h)\Big\}\, 
\frac{L(h+\frac{t}{\sigma_h})}{L(h)}
$$
and the log-Laplace transform
\be
\hat K_h(t) = K\Big(h+\frac{t}{\sigma_h}\Big) - K(h) - \frac{t}{\sigma_h}\,K'(h).
\en

In order to estimate (5.5) from above, assume that 
$K(h) \leq \frac{1}{2}\,h^2$, i.e. $A(h) \geq 0$ for all $h$. For $h \in A_n(a)$, 
the definition (4.5) implies that
$$
K(h) \geq \frac{1}{2}\,h^2 - \frac{a}{n-1},
$$
and hence
\begin{eqnarray}
\hat K_h(t) 
 & \leq &
\frac{1}{2}\, \big(h+t \sigma_h^{-1}\big)^2 - \frac{1}{2}\,h^2 + 
\frac{a}{n-1} - \frac{t}{\sigma_h}\,K'(h) \nonumber \\
 & = &
\frac 1 2 (t \sigma_h^{-1})^2 + \frac{a}{n-1} + t \sigma_h^{-1}\,(h - K'(h)).
\end{eqnarray}
Here the term $h - K'(h) = A'(h)$ can be estimated by virtue of the inequality (3.9),
which gives
$$
|h - K'(h)|^2 \leq 2 A(h) \leq \frac{2a}{n-1}
$$
and
$$
|t|\, \sigma_h^{-1}\,|h - K'(h)| \leq 
\frac{1}{2}\,(t \sigma_h^{-1})^2 + \frac{1}{2}\,|h - K'(h)|^2 \leq 
\frac{1}{2}\,(t \sigma_h^{-1})^2 + \frac{a}{n-1}.
$$
It follows from (5.6) that
$$
\hat K_h(t) \leq \frac{3}{2}\,(t \sigma_h^{-1})^2 + \frac{2a}{n-1}.
$$
Here, the right-hand side is bounded for sufficiently small $|t|$ and sufficiently
large $n$. One may require, for example, that $n \geq 4a + 1$ and
$|t| \leq \frac{1}{2} \sigma_h$, in which case
$\hat K_h(t) \leq 1$, so that
$$
\E\,e^{|t| \hat X(h)} \leq \E\,e^{t \hat X(h)} + \E\,e^{-t \hat X(h)} \leq 2 e.
$$
Using $x^3 e^{-|t| x} \leq (\frac{3}{e})^3\,|t|^{-3}$ ($x \geq 0$),
this gives
$\E\,|\widehat X(h)|^3 \leq 2e\,(\frac{3}{e})^3\,|t|^{-3}$.
One can summarize.

\vskip5mm
{\bf Lemma 5.2.} {\sl If the Laplace transform of  a subgaussian random 
variable $X$ is such that $A(h) \geq 0$ for all $h \in \R$, then for all $h \in A_n(a)$
with $n \geq 4a + 1$, we have
$$
\E\,e^{\sigma_h |\hat X(h)|/2} < 2e.
$$
As a consequence, $\E\,|\hat X(h)|^3 \leq C \sigma_h^{-3}$
up to some absolute constant $C>0$.
}

\vskip10mm
\section{{\bf Local Limit Theorem for Bounded Densities}}
\setcounter{equation}{0}

Before we can apply the representation (3.7), in the next step we need
to establish a uniform local limit theorem with a quantitative error term.
Let $(X_k)_{k \geq 1}$ be independent copies of a random variable $X$ 
with $\E X = 0$, $\Var(X) = 1$, $\beta_3 = \E\,|X|^3 < \infty$, which 
has a bounded density. Then the normalized sums $Z_n$
have bounded continuous densities $p_n$ for all $n \geq 2$ satisfying
$$
\sup_x |p_n(x) - \varphi(x)| = O\Big(\frac{1}{\sqrt{n}}\Big) \quad (n \rightarrow \infty).
$$
See for example \cite{Pe1,Pe2}. Let us quantify the error $O$-term 
in terms of $\beta_3$ and the maximum of density $M = M(X)$.

\vskip5mm
{\bf Lemma 6.1.} {\sl With some positive absolute constant $C$, we have
\be
\sup_x |p_n(x) - \varphi(x)| \leq C\,\frac{M^2 \beta_3}{\sqrt{n}}.
\en
}

{\bf Proof.} Denote by $f(t)$  the characteristic function of $X$. 
By the boundedness assumption, the characteristic functions
$$
f_n(t) = \E\,e^{itZ_n} = f(t/\sqrt{n})^n, \quad t \in \R,
$$
are integrable for all $n \geq 2$. Indeed, by the Plancherel theorem, 
$$
\int_{-\infty}^\infty |f(t)|^n\,dt \leq 
\int_{-\infty}^\infty |f(t)|^2\,dt = 2\pi \int_{-\infty}^\infty p(x)^2\,dx
\leq 2\pi M.
$$
Hence, one may apply the Fourier inversion
formula to represent the densities of $Z_n$ as
$$
p_n(x) = \frac{1}{2\pi} \int_{-\infty}^\infty e^{-itx} f_n(t)\,dt, \quad
x \in \R.
$$
Using a similar representation for the normal density, we get
$$
|p_n(x) - \varphi(x)| \leq 
\frac{1}{2\pi} \int_{-\infty}^\infty |f_n(t) - e^{t^2/2}|\,dt.
$$

As is well known (cf. e.g. \cite{Pe2}, p. 109),
$$
|f_n(t) - e^{t^2/2}| \leq 16\,\frac{\beta_3}{\sqrt{n}}\,|t|^3\,e^{-t^2/3}, \quad
|t| \leq \frac{\sqrt{n}}{4\beta_3},
$$
which yields
$$
\int_{|t| \leq \frac{\sqrt{n}}{4\beta_3}} |f_n(t) - e^{t^2/2}|\,dt \leq 
\frac{C\beta_3}{\sqrt{n}}
$$
with some absolute constant $C$. As for large values of $|t|$, it was shown 
in \cite{B-C-G4}, p.\,145,  that, for any $\ep \in (0,1]$ and $n \geq 4$
(which may be assumed in~(6.1)),
$$
\int_{|t| \geq \ep} |f(t)|^n\,dt \leq 
\frac{4\pi M}{\sqrt{2n}}\,\exp\big\{-n\ep^2/(5200 M^2)\big\}.
$$
Since $\beta_3 \geq 1$, this gives
$$
\int_{|t| \geq \frac{\sqrt{n}}{4\beta_3}} |f_n(t)|\,dt = \sqrt{n}
\int_{|t| \geq \frac{1}{4\beta_3}} |f(t)|^n\,dt \leq 
\frac{4\pi M}{\sqrt{2}}\,\exp\big\{-c_0 n/(\beta_3^2 M^2)\big\}.
$$
Since in general $M \geq 1/\sqrt{12}$, a similar estimate holds true for 
the normal density as well. As a result, we arrive at
$$
|p_n(x) - \varphi(x)| \leq C_0\,\Big(\frac{\beta_3}{\sqrt{n}} +
M \exp\{-c_0 n/(\beta_3^2 M^2)\}\Big)
$$
with some positive absolute constants $C_0$ and $c_0$,
Using $e^{-x^2} < x^{-1}$ ($x>0$), the second term in
the brackets is dominated by the first one up to the multiple of $M^2$.
Hence, the above estimate may be simplified to (6.1).
\qed

\section{{\bf Local Limit Theorem for Shifted Densities}}
\setcounter{equation}{0}

An application of Lemma 6.1 to the normalized sums of independent
copies of random variables $\hat X(h)$ defined in (5.4) leads to
the following refinement of the representation (3.7) from Proposition 3.3,
when the point $x$ belongs to the critical zone $A(x) \leq \frac{a}{n-1}$. Define
$$
v_x = \frac{x - m_x}{\sigma_x} = \frac{x - K'(x)}{\sigma_x} = \frac{A'(x)}{\sigma_x},
$$
where we recall that $m_x = K'(x)$ and $\sigma_x^2 = K''(x)$.

\vskip5mm
{\bf Lemma 7.1.} {\sl If the Laplace
transform of a subgaussian random variable $X$ with finite constant
$c = 1 + T_\infty(p||\varphi)$
is such that $A(h) \geq 0$ for all $h \in \R$, then for all $x \in A_n(a)$
with $n \geq 4(a + 1)$, we have
\be
\frac{p_n(x\sqrt{n})}{\varphi(x\sqrt{n})} = 
\frac{1}{\sigma_x}\,e^{-nA(x)-nv_x^2/2} + \frac{Bc^4}{\sqrt{n}},
\en
where $B = B_n(x)$ is bounded by an absolute constant.
}

\vskip5mm
{\bf Proof.} Let us return to the term $Q_{h_n} p_n$ in (3.7) with
$h_n = h\sqrt{n}$. By Proposition 3.2, this density has a convolution structure.
Recall that, for any random variable $X$ with density $p = p_X$,
$$
Q_h p_{\lambda X}(x) = \frac{1}{\lambda}\,(Q_{\lambda h} p)\Big(\frac{x}{\lambda}\Big).
$$
Using this notation, $p_n = p_{S_n/\sqrt{n}}$ in terms of the sum
$S_n = X_1 + \dots + X_n$. Hence with $\lambda = 1/\sqrt{n}$, 
$$
Q_{h_n} p_n(x) = \sqrt{n}\, (Q_h p_{S_n})(x\sqrt{n}) = 
\sqrt{n}\, (Q_h p) * \dots * (Q_h p) (x\sqrt{n}),
$$
where we applied Proposition 3.2 in the last step. By definition, $Q_h p$ is
the density of the random variable $X(h)$. Hence, $Q_{h_n} p_n(x)$ represents
the density for the normalized sum
$$
Z_{n,h} \equiv (X_1(h) + \dots + X_n(h))/\sqrt{n},
$$
assuming that $X_k(h)$ are independent. Introduce the normalized sums
\be
\hat Z_{n,h} \equiv (\hat X_1(h) + \dots + \hat X_n(h))/\sqrt{n}
\en
for the shifted distributions (5.4), i.e. with
$X_k(h) = m_h + \sigma_h \hat X_k(h)$. Thus,
$$
Z_{n,h} = m_h \sqrt{n} + \sigma_h \hat Z_{n,h}.
$$

Denote by $\hat p_{n,h}$ the density of $\hat Z_{n,h}$. Then the density of
$Z_{n,h}$ is given by
$$
p_{n,h}(x) = 
\frac{1}{\sigma_h}\,\hat p_{n,h}\Big(\frac{x - m_h \sqrt{n}}{\sigma_h}\Big),
\quad x \in \R.
$$
At the points $x_n = x\sqrt{n}$ as in (3.7), we therefore obtain that
$$
Q_{h_n} p_n(x_n) = p_{n,h}(x_n) = 
\frac{1}{\sigma_h}\,\hat p_{n,h}\Big(\frac{x - m_h}{\sigma_h} \sqrt{n}\Big).
$$
Consequently, the equality (3.7) may be equivalently stated as
$$
\frac{p_n(x\sqrt{n})}{\varphi(x\sqrt{n})} = 
\sqrt{2\pi}\,e^{\frac{n}{2}\, (x - h)^2 - n A(h)}\,
\frac{1}{\sigma_h}\,\hat p_{n,h}\Big(\frac{x - m_h}{\sigma_h} \sqrt{n}\Big).
$$
In particular, for $h = x$, we get
\be
\frac{p_n(x\sqrt{n})}{\varphi(x\sqrt{n})} = \sqrt{2\pi}\,e^{- n A(x)}\,
\frac{1}{\sigma_x}\,\hat p_{n,x}(v_x\sqrt{n}).
\en

We are now in a position to apply Lemma 6.1 to the sequence $\hat X_k(x)$ and write
\be
\hat p_{n,x}(z) = \varphi(z) + B\,\frac{\beta_3(x) M(x)^2}{\sqrt{n}},
\quad z \in \R,
\en
where the quantity $B = B_n(z)$ is bounded by an absolute constant, 
$\beta_3(x) = \E\,|\hat X(x)|^3$, and $M(x) = M(\hat X(x))$.  The latter
maximum can be bounded by virtue of the upper bound (4.2):
$$
M(\hat X(x)) = \sigma_x M(X(x)) = \sigma_x M(Q_x p)
\leq \frac{c\sigma_x}{\sqrt{2\pi}}\, e^{A(x)}.
$$
In this case, (7.4) may be simplified with a new $B$ to
$$
\hat p_{n,x}(z) = \varphi(z) + Bc^2\,\frac{\beta_3(x) \sigma_x^2}{\sqrt{n}}\, e^{2A(x)}.
$$
Inserting this in (7.3) with $z = v_x \sqrt{n}$, we arrive at
$$
\frac{p_n(x\sqrt{n})}{\varphi(x\sqrt{n})} = \frac{1}{\sigma_x}\,
e^{-nA(x) -nv_x^2/2} + Bc^2\,\frac{\beta_3(x) \sigma_x}{\sqrt{n}}\, e^{-(n-2)A(x)}.
$$

To further simplify, assume that $x \in A_n(a)$ with $n \geq 4(a + 1)$. Then, 
by Lemmas 5.1-5.2, $\beta_3(x) \leq C\sigma_x^{-3}$, while 
$\sigma_x^{-1} \leq 2c\, e^{A(x)}$. Hence,
$$
\beta_3(x)\sigma_x\, e^{-(n-2)A(x)} \leq 4C c^2\, e^{-(n-4)A(x)} \leq 4C c^2.
$$
\qed

\section{{\bf Proof of Theorem 2.1}}
\setcounter{equation}{0}

Recall that the assumptions 1)-2) stated before Theorem 2.1 are necessary
for the convergence $T_\infty(p_n||\varphi) \rightarrow 0$ as $n \rightarrow \infty$.
For simplicity, we assume that $n_0 = 1$, that is, $X$ is a strictly subgaussian random 
variable with mean zero, variance one, and with finite constant $c = 1 + T_\infty(p||\varphi)$. 
In particular, the function
$$
A(x) = \frac{1}{2} x^2 - K(x)
$$
is non-negative on the whole real line. 

{\bf Sufficiency part.} The critical zones $A_n(a) = \{x \in \R: A(x) \leq \frac{a}{n-1}\}$ was
defined for a parameter $a>0$ and $n \geq 2$. Choosing $a = \log(1/\ep)$ for a given $\ep \in (0,1)$,
we have, by (4.6),
\be
\sup_{x \notin A_n(a)}\, \frac{p_n(x\sqrt{n})}{\varphi(x\sqrt{n})} \leq c\sqrt{2}\,\ep.
\en

In the case $x \in A_n(a)$ with $n \geq 4(a + 1)$, the equality (7.1) is applicable
and implies
$$
\sup_{x \in A_n(a)}\, \frac{p_n(x\sqrt{n})}{\varphi(x\sqrt{n})} \leq 
\sup_{x \in A_n(a)}\, \frac{1}{\sigma_x} + O\Big(\frac{1}{\sqrt{n}}\Big).
$$
Using (8.1), we conclude that
$$
1 + T_\infty(p_n||\varphi) \leq \sup_{x \in A_n(a)}\, \frac{1}{\sigma_x} + 
c\sqrt{2}\,\ep + O\Big(\frac{1}{\sqrt{n}}\Big).
$$
Thus, a sufficient condition for the convergence $T_\infty(p_n||\varphi) \rightarrow 0$ 
as $n \rightarrow \infty$ is that, for any $\ep \in (0,1)$,
$$
\limsup_{n \rightarrow \infty}\, \sup_{x \in A_n(\log(1/\ep))}\, 
\sigma_x^{-2} \leq 1.
$$
Equivalently, we need to require that
$\liminf_{n \rightarrow \infty}\, \inf_{x \in A_n(a)}\, K''(x) \geq 1$
for any $a>0$, that is,
$$
\limsup_{n \rightarrow \infty}\, \sup_{x \in A_n(a)}\, A''(x) \leq 0.
$$
Since $A(x) = O(1/n)$ on every set $A_n(a)$, the above may be written as
the following continuity condition
\be
\lim_{A(x) \rightarrow 0} \max(A''(x),0) = 0.
\en

{\bf Necessity part.} 
To see that the condition (8.2) is also necessary for the convergence in $T_\infty$,
let us return to the representation (7.1). Assuming that 
$T_\infty(p_n||\varphi) \rightarrow 0$, it implies that, for any $a>0$,
\be
\limsup_{n \rightarrow \infty}\, \sup_{x \in A_n(a)} \frac{1}{\sigma_x}
\exp\Big\{- n \Big(A(x) + \frac{1}{2}v_x^2\Big)\Big\} \leq 1.
\en
Recall that
$$
A'(x)^2 \leq 2A(x), \quad \sigma_x^{-2} \leq \frac{6}{\pi}\,c^2 e^{A(x)}.
$$
(cf. Remark 3.5 and Lemma 5.1). Hence
$$
v_x^2 = \frac{A'(x)^2}{\sigma_x^2} \leq \frac{2A(x)}{\sigma_x^2} \leq
\frac{12}{\pi}\,c^2 e^{A(x)} A(x) \leq 12\,c^2 A(x),
$$
assuming that $x \in A_n(a)$ with $a \leq 1$ and $n \geq 2$ 
in the last step. Since  $nA(x) \leq 2a$ on the set $A_n(a)$ and $c \geq 1$,
it follows that 
$$
A(x) + \frac{1}{2}v_x^2 \leq 7 c^2 A(x) \leq \frac{14 c^2}{n}\,a.
$$
Thus, (8.3) implies that
$$
\limsup_{n \rightarrow \infty}\, \sup_{x \in A_n(a)} \frac{1}{\sigma_x}
\leq e^{14 c^2 a}, \quad 0 < a \leq 1.
$$
Therefore, for all $n \geq n(a)$,
$$
\inf_{x \in A_n(a)} K''(x) \geq e^{-30 c^2 a}.
$$
Since $a$ may be as small as we wish, we conclude that, for any $\ep > 0$, there is 
$\delta>0$ such that $A(x) \leq \delta \, \Rightarrow \, K''(x) \geq 1-\ep$, or
$A(x) \leq \delta \, \Rightarrow \, A''(x) \leq \ep$. But this is the same as (8.2).
\qed

\vskip2mm
One wide class of strictly subgaussian distributions with mean zero and variance one
is described in terms of the Laplace transform $L(t) = \E e^{tX}$
via the potential requirement (2.2), i.e.
\be
L(t) \leq (1 - \delta)\, e^{t^2/2}
\en
for all $t_0 > 0$ and $|t| \geq t_0$ with some $\delta = \delta(t_0)$, $\delta \in (0,1)$.
In this case, the log-Laplace transform and the $A$-function satisfy
$$
K(t) \leq \frac{1}{2} t^2 + \log(1 - \delta), \quad A(t) \geq - \log(1 - \delta).
$$
Hence, the approach $A(t) \rightarrow 0$ is only possible when $t \rightarrow 0$.
But, for strictly subgaussian distributions, we necessarily have $A(t) = O(t^4)$ and $A''(t) = O(t^2)$
near zero. Therefore, the condition (8.2) is fulfilled automatically.

\vskip5mm
{\bf Corollary 8.1.} {\sl If a random variable $X$ with mean zero, variance one,
and finite distance $T_\infty(p||\varphi)$ satisfies the separation property $(8.4)$,
then $T_\infty(p_n||\varphi) \rightarrow 0$ as $n \rightarrow~\infty$.
}

\vskip5mm
\section{{\bf Characterization in the Periodic Case. Examples}}
\setcounter{equation}{0}

Let us apply Theorem 2.1 to the Laplace transforms $L(t)$
with $L(t) e^{-t^2/2}$ being periodic.
Suppose that $\E X = 0$ and $\Var(X) = 1$. As before, assume that:

\vskip2mm
1) $Z_n$ has density $p_n$ for some $n=n_0$ such that 
$T_\infty(p_n||\varphi) < \infty$;

2) $X$ is strictly subgaussian: $L(t) \leq e^{t^2/2}$,
or equivalently $\Psi(t) \leq 1$ for all $t \in \R$, where 
\be
\Psi(t) = L(t)\,e^{-t^2/2}, \quad t \in \R.
\en

\noindent
In addition, suppose that the function $\Psi(t)$ is $h$-periodic for some $h>0$.

\vskip5mm
{\bf Proof of Theorem 2.3} (first part). We need to show that the convergence
$T_\infty(p_n||\varphi) \rightarrow 0$ is equivalent
to the assertion that, for every $0 < t < h$,
\be
\Psi(t) = 1 \, \Rightarrow \, \Psi''(t) = 0.
\en

First note that, due to $\Psi(t)$ being analytic, the equation
$\Psi(t) = 1$ has finitely many solutions in the interval $[0,h]$ only,
including the points $t = 0$ and $t = h$ (by the periodicity).
Hence, the condition $b)$ in Theorem 2.1 may be ignored, and we obtain
that $T_\infty(p_n||\varphi) \rightarrow 0$ as $n \rightarrow \infty$,
if and only if 
\be
A''(t) = 0 \ {\rm for \ every \ point} \ t \in [0,h] \ {\rm such \ that} \ A(t) = 0.
\en
Here one may exclude the endpoints, since $A''(0) = A''(h) = 0$, by the strict
subgaussianity and periodicity. As for the interior points $t \in (0,h)$, note that 
$A(t) = -\log \Psi(t)$ has the second derivative
$$
A''(t) = \frac{\Psi'(t)^2 - \Psi''(t) \Psi(t)}{\Psi(t)^2} = -\Psi''(t)
$$
at every point $t$ such that $\Psi(t) = 1$ (in which case $\Psi'(t) = 0$
due to the property $\Psi \leq 1$). This shows that (9.3) is reduced to the
condition (9.2).
\qed

\vskip2mm
In order to describe examples illustrating Theorem 2.3, let us start with the following.

\vskip2mm
{\bf Definition.} We say that the distribution $\mu$ of a random 
variable $X$ is periodic with respect to the standard normal law, 
with period $h>0$, if it has a density $p(x)$ such that the function
$$
q(x) = \frac{p(x)}{\varphi(x)} = \frac{d\mu(x)}{d\gamma(x)},
\quad x \in \R,
$$
is periodic with period $h$, that is, $q(x+h) = q(x)$ for all $x \in \R$.

\vskip5mm
Here, $q$ represents the density of $\mu$ with respect to the standard Gaussian 
measure $\gamma$. We denote the class of all such distributions by 
$\mathfrak F_h$, and say that $X$ belongs to $\mathfrak F_h$.
Let us briefly collect and recall without proof several observations from 
\cite{B-C-G5} on this interesting class of probability distributions
(cf. Sections 10-13). 

\vskip5mm
{\bf Proposition 9.1.} {\sl If $X$ belongs to $\mathfrak F_h$, then $X$ is 
subgaussian, and the function $\Psi(t)$ in $(9.1)$ is $h$-periodic.
It may be extended to the complex plane as an entire function. 
Conversely, if $\Psi(t)$ for a subgaussian random variable
$X$ is $h$-periodic, then $X$ belongs to $\mathfrak F_h$, as long as 
the characteristic function $f(t)$ of $X$ is integrable. 
}

\vskip5mm
Since
$$
f(t) = L(it) = \Psi(it)\,e^{-t^2/2},
$$
the integrability assumption in the reverse statement is fulfilled,
as long as $\Psi(z)$ has order smaller than 2, that is, when
$|\Psi(z)| = O(\exp\{|z|^\rho\})$ as $|z| \rightarrow \infty$ for some $\rho<2$.

The periodicity property is stable along convolution: The normalized sums $Z_n$
belong to $\mathfrak F_{h\sqrt{n}}$, as long as $X$ belongs to $\mathfrak F_h$.

This class contains distributions whose Laplace transform has the form 
$L(t) = \Psi(t)\, e^{t^2/2}$, where $\Psi$ is a trigonometric polynomial. 
More precisely, consider functions of the form
$$
\Psi(t) = 1 - c P(t), \quad P(t) = a_0 + \sum_{k=1}^N (a_k \cos(kt) + b_k \sin(kt)), 
$$
where $a_k,b_k$ are given real coefficients, and $c \in \R$
is a non-zero parameter. 

\vskip5mm
{\bf Proposition 9.2.} {\sl If $P(0) = 0$ and $|c|$ is small enough, then 
$L(t)$ represents the Laplace transform of a subgaussian random variable $X$
with density $p(x) = q(x) \varphi(x)$, where $q(x)$ is
a non-negative trigonometric polynomial of degree at most $N$.
}

\vskip5mm
Note that necessarily $q$ is bounded, so that $T_\infty(p||\varphi) < \infty$.
As for the requirement that $P(0) = a_0 + \sum_{k=1}^N a_k = 0$,
it guarantees that $\int_{-\infty}^\infty p(x)\,dx = 1$. 
In order to apply Theorem 2.3, there are two more constraints
coming from the assumption that $\E X = 0$ and $\E X^2 = 1$.

\vskip5mm
{\bf Corollary 9.3.} {\sl Suppose that the polynomial $P(t)$ satisfies

\vskip2mm
$1)$ $P(0) = P'(0) = P''(0) = 0$;

$2)$ $P(t) \geq 0$ for $0<t<h$, where $h$ is the smallest period of $P$.

\vskip2mm
\noindent
If $c>0$ is small enough, then $L(t)$ represents the Laplace transform 
of a strictly subgaussian random variable $X$. Moreover, if $P(t) > 0$ for 
$0<t<h$, then $T_\infty(p_n||\varphi) \rightarrow 0$ as $n \rightarrow \infty$.
}

\vskip2mm
In terms of the coefficients of the polynomial, the moment assumptions
$P'(0) = P''(0) = 0$ are equivalent to
$\sum_{k=1}^N k b_k = \sum_{k=1}^N k^2 a_k = 0$.
The assumption 2) implies that $0 < \Psi(t) \leq 1$, and if $P(t) > 0$ for 
$0<t<h$, then the equation $\Psi(t)=1$ has no solution in this interval.

\vskip5mm
{\bf Example 9.4.} Consider the transforms of the form
\be
L(t) = (1 - c \sin^m(t))\, e^{t^2/2}
\en
with an arbitrary integer $m \geq 3$, where $|c|$ is small enough.
Then $\E X = 0$, $\E X^2 = 1$, and the cumulants of $X$ satisfy
$\gamma_k(X) = 0$ for all $3 \leq k \leq m-1$.

Moreover, if $m \geq 4$ is even, and $c>0$, the random variable 
$X$ with the Laplace transform (9.4) is strictly subgaussian.
Hence the conditions in Corollary 9.3 are met, and we obtain
the statement about the R\'enyi divergence of infinite order.
In the case $m=4$, (9.4) corresponds to
$$
P(t) = \sin^4 t = \frac{1}{8}\,(3 - 4 \cos(2t) + \cos(4t)).
$$

{\bf Example 9.5.} Put
\be
P(t) = (1 - 4\sin^2 t)^2 \, \sin^4 t.
\en
Then, $P(t) = O(t^4)$, implying that
$P(0) = P'(0) = P''(0) = 0$. Note that $\Psi(t) = 1 - cP(t)$
is $\pi$-periodic, and $h = \pi$ is the smallest period, although
$$
\Psi(0) = \Psi(t_0) = \Psi(\pi) = 1, \quad t_0 = \pi/6.
$$
As we know, if $c>0$ is small enough, then $L(t) = 1 - c\Psi(t)$
represents the Laplace transform of a strictly subgaussian random 
variable $X$. In this case, the last assertion in Corollary 9.3 is not applicable.
Thus, the property that $h$ is the smallest period for a periodic 
function $\Psi(t)$ such that $0 \leq \Psi(t) \leq 1$ and 
$\Psi(h) = 1$ does not guarantee that $0<\Psi(t) < 1$ for $0 < t < h$.

Nevertheless, all assumptions of Theorem 2.3 are fulfilled for sufficiently 
small $c>0$ with $h = \pi$, and we may check the condition (9.2). In this case, 
$$
\Psi(t) = 1 - c Q(t)^2, \quad Q(t) = (1 - 4\sin^2 t) \, \sin^2 t =  \sin^2 t - 4\sin^4 t,
$$
so that
$$
\Psi''(t) = -2c\,(Q(t) Q''(t) + Q'(t)^2) = -2c Q'(t)^2
$$
at the points $t$ such that $Q(t) = 0$, that is, for $t = t_0$. Hence
$\Psi''(t) = 0 \Leftrightarrow Q'(t) = 0$. In our case,
$$
Q'(t) = 2 \sin t \cos t - 16 \sin^3 t \cos t = \sin(2t)\,(1 - 8 \sin^2 t),
$$
$$
Q'(t_0) = \sin(\pi/3)\,(1 - 8 \sin^2(\pi/6)) = -\frac{\sqrt{3}}{2} \neq 0.
$$
Hence $\Psi''(t_0) \neq 0$, showing that the condition (9.2) is {\bf not} fulfilled.
Thus, the CLT with respect to $T_\infty$ does not hold in this example.

The examples based on trigonometric polynomials may be generalized to the setting
of $2\pi$-periodic functions represented by Fourier series
$$
P(t) = a_0 + \sum_{k=1}^\infty (a_k \cos(kt) + b_k \sin(kt)). 
$$
Then, the assertions in Proposition 9.2 and Corollary 9.3 will continue to hold,
as long as the coefficients satisfy 
$\sum_{k=1}^\infty e^{k^2/2} (|a_k| + |b_k|) < \infty$.

\vskip10mm
\section{{\bf Richter's Local Limit Theorem and its Refinement}}
\setcounter{equation}{0}

We now turn to the problem of convergence rates with respect to $T_\infty$,
which can be explored, for example, under the separation-type condition (2.2).
In this case, it was shown in Corollary 4.3 that $p_n(x)$ is much smaller than 
$\varphi(x)$ outside the interval $|x| = O(\sqrt{n})$. In the region 
$|x| = o(\sqrt{n})$, an asymptotic behavior of the densities $p_n$ 
of the normalized sums 
$$
Z_n = (X_1 + \dots + X_n)/\sqrt{n}
$$
is governed by the following theorem due to Richter \cite{Ri}. 
Assume that $(X_n)_{n \geq 1}$ are independent copies of a random
variable $X$ with mean $\E X = 0$ and variance $\Var(X) =~1$.

\vskip5mm
{\bf Theorem 10.1.} {\sl Let $\E\,e^{c|X|} < \infty$ for some 
$c>0$, and let $Z_n$ have a bounded density for some 
$n$. Then $Z_n$ with large $n$ have bounded continuous
densities $p_n$ satisfying
\be
\frac{p_n(x)}{\varphi(x)} = \exp\Big\{\frac{x^3}{\sqrt{n}}\,
\lambda\Big(\frac{x}{\sqrt{n}}\Big)\Big\}\,
\Big(1 + O\Big(\frac{1+|x|}{\sqrt{n}}\Big)\Big)
\en
uniformly for $|x| = o(\sqrt{n})$. The function
$\lambda(z)$ is represented by an infinite power series
which is absolutely convergent in a neighborhood of $z=0$.
}

\vskip5mm
The corresponding representation
\be
\lambda(z) = \sum_{k=0}^\infty \lambda_k z^k
\en
is called Cramer's series; it is analytic in some disc $|z| \leq \tau_0$ 
of the complex plane. The proof of this theorem may also be found 
in the book by Ibragimov and Linnik \cite{I-L}, cf. Theorem 7.1.1, 
where it was assumed that $X$ has a continuous bounded density. 
The representation (10.1) was further investigated 
there for zones of normal attraction
$|x| = o(n^\alpha)$, $\alpha < \frac{1}{2}$.

One immediate consequence of (10.1) is that
\be
\frac{p_n(x)}{\varphi(x)} \rightarrow 1 \quad 
{\rm as} \ n \rightarrow \infty
\en
uniformly in the region $|x| = o(n^{1/6})$. However,
in general this is no longer true outside this region.
To better understand the possible behavior of densities,
one needs to involve the information about the coefficients 
in the power series (10.2). As was already mentioned in~\cite{I-L},
$\lambda_0 = \frac{1}{6}\,\gamma_3$,
$\lambda_1 = \frac{1}{24}\,(\gamma_4 - 3\gamma_3^2)$.
However, in order to judge the behavior $\lambda(z)$ for small
$z$, one should describe the leading term in this series.
The analysis of the saddle point associated to the
log-Laplace transform of the distribution of $X$ shows that
\be
\lambda(z) = \frac{\gamma_m}{m!}\,z^{m-3} + O(|z|^{m-2}), \quad
{\rm as} \ z \rightarrow 0,
\en
where $\gamma_m$ denotes the first non-zero cumulant
of $X$ (when $X$ is not normal). Equivalently, $m$ is 
the smallest integer such that $m \geq 3$ and
$\E X^m \neq \E Z^m$, where $Z \sim N(0,1)$. In this case
$\gamma_m = \E X^m - \E Z^m$.

Using (10.4) in (10.1), we obtain a more informative representation
\be
\frac{p_n(x)}{\varphi(x)} = \exp\bigg\{\frac{\gamma_m}{m!}\, 
\frac{x^m}{n^{\frac{m}{2} - 1}} + 
O\Big(\frac{x^{m+1}}{n^{\frac{m}{2}}}\Big)\bigg\}
\Big(1 + O\Big(\frac{1+|x|}{\sqrt{n}}\Big)\Big),
\en
which holds uniformly for $|x| = o(\sqrt{n})$. With this 
refinement, the convergence in (10.3) holds true uniformly over all $x$ 
in the potentially larger region
$$
|x| \leq \ep_n\,n^{\frac{1}{2} - \frac{1}{m}} \quad
(\ep_n \rightarrow 0).
$$ 
For example, if the distribution of $X$ is symmetric about the origin,
then $\gamma_3 = 0$, so that necessarily $m \geq 4$.

Nevertheless, for an application to the $T_\infty$-distance, it is desirable
to get some information for larger intervals such as 
$|x| \leq \tau_0\sqrt{n}$ and to replace the term 
$O(\frac{|x|}{\sqrt{n}})$ in (10.5) with an explicit $n$-dependent quantity. 
To this aim, the following refinement of Theorem 10.1 was 
recently proved in \cite{B-C-G4}.

\vskip5mm
{\bf Theorem 10.2.} {\sl Let the conditions of Theorem $10.1$
be fulfilled. There is $\tau_0 > 0$ with the following property.
Putting $\tau = x/\sqrt{n}$, we have for $|\tau| \leq \tau_0$
\be
\frac{p_n(x)}{\varphi(x)} = e^{n\tau^3 \lambda(\tau) - \mu(\tau)}\,
\big(1 + O(n^{-1} (\log n)^3)\big),
\en
where $\mu(\tau)$ is an analytic function in $|\tau| \leq \tau_0$
such that $\mu(0)=0$.
}

\vskip2mm
Here, similarly to (10.4),
$$
\mu(\tau) = 
\frac{1}{2(m-2)!}\,\gamma_m \tau^{m-2} + O(|\tau|^{m-1}).
$$

As a consequence of (10.6), which cannot be obtained on the basis
of (10.1) or (10.5), we have the following assertion which was also derived 
in \cite{B-C-G4}.

\vskip5mm
{\bf Corollary 10.3.} {\sl Under the same conditions, suppose that the first non-zero cumulant $\gamma_m$ of $X$ is negative and
$m$ is even. There exist constants 
$\tau_0 > 0$ and $c>0$ with the following property.
If $|\tau| \leq \tau_0$, $\tau = x/\sqrt{n}$, then
\be
\frac{p_n(x)}{\varphi(x)} \leq 1 + \frac{c(\log n)^3}{n}.
\en
}

{\bf Proof of Theorem 2.2.} It remains to combine Corollary 4.3 with 
Corollary 10.3 and note that, for any strictly subgaussian random variable 
$X$ with variance one, $m$ is even and $\gamma_m < 0$. 
Indeed, the log-Laplace transform of the distribution of $X$ admits 
the following Taylor expansion near zero
$$
K(t) = \log \E\,e^{tX} = \frac{1}{2}\,t^2 + 
\sum_{k=3}^\infty \frac{\gamma_k}{k!}\,t^k = \frac{1}{2}\,t^2 +
\frac{\gamma_m}{m!}\,t^m + O(t^{m+1}),
$$ 
which is a definition of cumulants.
Hence, the strict subgaussianity, that is, the property $K(t) \leq \frac{1}{2}\,t^2$ 
for all $t \in \R$ implies the claim.
\qed

\vskip5mm
{\bf Proof of Theorem 2.3} (convergence part). 
For simplicity, let $n_0 = 1$, so that the random variable $X$ has density $p$ 
with $T_\infty(p||\varphi) < \infty$. By the assumption, $\E X = 0$, $\Var(X) = 1$,
and
$$
L(t) = \E\,e^{tX} = \Psi(t)\,e^{t^2/2}, \quad t \in \R,
$$
for some periodic function $\Psi(t)$ with period $h>0$ such that
$0 < \Psi(t) < 1$ for all $0 < t < h$. Hence
$$
L(t/\sqrt{n})^n = \E\,e^{tZ_n} = \Psi_n(t)\,e^{t^2/2}, \quad
\Psi_n(t) = \Psi(t/\sqrt{n})^n,
$$
where the function $\Psi_n(t)$ has period $h\sqrt{n}$. Since the density
$p$ is bounded, the characteristic function of $X$ is square
integrable. Hence, the characteristic function of $Z_n$ is integrable
whenever $n \geq 2$. In this case, we are in position to apply Proposition 9.1 
to the random variable $Z_n$ and conclude that it has a continuous density $p_n$ 
which is periodic with respect to the standard normal law with period $h\sqrt{n}$. 
That is, $p_n(x) = q_n(x) \varphi(x)$
for some continuous, periodic function $q_n$ with period $h\sqrt{n}$.
We need to show that 
\be
\sup_x \, (q_n(x) - 1) = O\Big(\frac{(\log n)^3}{n}\Big) \quad
{\rm as} \ n \rightarrow \infty.
\en

In view of periodicity, one may restrict this supremum to the interval 
$0 \leq x \leq h\sqrt{n}$. But, if $0 \leq x \leq \tau_0\sqrt{n}$, where 
$\tau_0$ is taken as in Corollary 10.3, we obtain the desired rate due to (10.7). 
Here, without loss of generality one may 
assume that $\tau_0 < h$. Since $q_n(x) = q_n(x - h\sqrt{n})$,
the same conclusion is also true, if we restrict the supremum 
to $(h - \tau_0)\sqrt{n} \leq x \leq h\sqrt{n}$.
Finally, if $\tau_0\sqrt{n} \leq x \leq (h - \tau_0)\sqrt{n}$,
we apply the bound (4.3) which gives
$$
q_n(x) \leq c\sqrt{2}\,\Psi\Big(\frac{x}{\sqrt{n}}\Big)^{n-1}, \quad
c = 1 + T_\infty(p||\varphi).
$$
Since $\Psi(t)$ is continuous, $\sup_{\tau_0 \leq t \leq h - \tau_0} \Psi(t) < 1$.
Hence the expression on the right-hand side is exponentially
small for growing $n$. Collecting these estimates, we arrive at (10.8).
\qed

\section{{\bf Examples Based on Weighted Sums}}
\setcounter{equation}{0}

Here we describe some examples illustrating Theorem 2.2. It
involves the separation condition (2.2) on the Laplace transform,
\be
\sup_{|t| \geq t_0} \big[e^{-t^2/2}\, \E\,e^{tX} \big] < 1 \quad 
{\rm for \ all} \ t_0>0,
\en
and
states the following  speed of convergence  in the CLT
\be
D_\infty(p_n||\varphi) = O\Big(\frac{(\log n)^3}{n}\Big)
\quad {\rm as} \ \ n \rightarrow \infty,
\en
provided that the  necessary condition $D_\infty(p_n||\varphi) < \infty$
for some $n = n_0$ holds, where
$p_n$ denote the densities of the normalized sums $Z_n$
constructed for independent copies of 
a random variable $X$ with $\E X = 0$, $\Var(X) = 1$.

While in general this condition is rather delicate,
in the simplest case $n_0 = 1$, it reduces to the pointwise subgaussian
bound
\be
p(x) \leq M \varphi(x), \quad x \in \R,
\en
which should hold with some constant $M$ for a density $p$ of the random 
variable $X$. This property is obviously fulfilled, when the density $p$ 
is bounded and compactly supported; the rate  (11.2) holds as well
for a family of probability distributions whose Laplace 
transform contains a periodic component (see remarks after Proposition 9.2). 
We now consider further examples 
where the density $p$ is representable as a ``weighted" convolution of 
at least two densities satisfying (11.3). More precisely, we have:

\vskip4mm
{\bf Corollary 11.1.} {\sl Assume that $X$
satisfies $(11.1)$ and is represented as
\be
X = c_0 \eta_0 + c_1 \eta_1 + c_2 \eta_2, \quad
c_0^2 + c_1^2 + c_2^2 = 1, \ \ c_1,c_2 > 0, 
\en
where the independent random variables $\eta_k, k=0,1,2$ are strictly subgaussian 
with variance one and satisfy $D_\infty(\eta_k||\varphi) < \infty$ for 
$k = 1,2$. Then the CLT holds with rate $(11.2)$.
}

\vskip4mm
As an interesting subclass, one may consider infinite
weighted convolutions, that is, random variables of the form
\be
X = \sum_{k=1}^\infty a_k \xi_k, \quad 
\sum_{k=1}^\infty a_k^2 = 1.
\en

{\bf Corollary 11.2.} {\sl Assume that the i.i.d. random variables 
$\xi_k$ are strictly subgaussian and have a bounded, compactly 
supported density with variance $\Var(\xi_1) = 1$. If $\xi_1$ 
satisfies $(11.1)$, then the CLT holds with rate $(11.2)$. 
}

\vskip4mm
This statement includes, for example, infinite weighted 
convolutions of the uniform distribution on a bounded symmetric interval.

By Theorem 2.2, Corollary 11.1 follows from the next
general assertion.

\vskip4mm
{\bf Lemma 11.3.} {\sl Suppose that the random variable $X$ is 
represented in the form $(11.4)$, where the random variables
$\eta_0,\eta_1,\eta_2$ are independent and possess the properties:

\vskip2mm
$a)$ $\eta_0$ is strictly subgaussian with $\Var(\eta_0)=1$;

$b)$ $\eta_1,\eta_2$ have densities $q_1,q_2$
such that $q_k(x) \leq M_k \varphi(x)$ for all $x \in \R$
with some constants $M_k$ $(k = 1,2)$.

\vskip2mm
\noindent
Then $X$ has a density $p$ satisfying $(11.3)$
with constant $M = \frac{1}{\sqrt{2c_1 c_2}}\,M_1 M_2$.
}

\vskip4mm
{\bf Proof.} The case $c_0 = 0$ is simple. Then $X$ has density
$$
p(x) = \frac{1}{c_1 c_2} \int_{-\infty}^\infty
q_1\Big(\frac{x-y}{c_1}\Big)\, q_2\Big(\frac{y}{c_1}\Big)\,dy, \quad
x \in \R,
$$
which, by the assumption, does not exceed
$$
\frac{M_1 M_2}{c_1 c_2} \int_{-\infty}^\infty
\varphi\Big(\frac{x-y}{c_1}\Big)\, \varphi\Big(\frac{y}{c_1}\Big)\,dy =
M_1 M_2\,\varphi(x).
$$
Hence, (11.3) is fulfilled with constant $M = M_1 M_2$
(which is better than what is claimed in the lemma, since $2c_1 c_2 \leq 1$).

In the basic case $c_0>0$, introduce the characteristic functions 
$f_k(t)$ of $\eta_k$ and put $g_k(t) = f_k(c_k t)$, $k = 0,1,2$. 
Since the densities $q_1,q_2$ are bounded, they belong to $L^2(\R)$ 
together with their characteristic functions $f_1,f_2$, according 
to the Plancherel theorem. The same is true for $g_1,g_2$, so that 
the characteristic function of $X$,
\be
f(t) = g_0(t) g_1(t) g_2(t),
\en
is integrable on the real line (using $|g_0(t)| \leq 1$ for all $t \in \R$). 
As a consequence, the random variable $X$ has a continuous density described 
by the inversion formula
\be
p(x) = \frac{1}{2\pi} \int_{-\infty}^\infty e^{-itx} f(t)\,dt, \quad
x \in \R.
\en

Moreover, the pointwise subgaussian bounds on the densities $q_k$ in 
$b)$ for $k=1,2$ ensure that $\E\,e^{\lambda \eta_k^2} < \infty$ for 
$\lambda < \frac{1}{2}$, implying that the random variables $\eta_k$ 
are subgaussian. Since $\eta_0$ is also subgaussian, we conclude that 
$X$ is subgaussian as well. Hence, all $g_k(t)$ and $f(t)$
may be extended 
from the real line to the complex plane as entire functions 
of order at most 2, and thus, (11.6) holds true for all $t \in \C$.

For definiteness, let $x < 0$ in (11.7). We use a contour integration
to obtain a different representation for $p(x)$.
Fix $T>0$, $y>0$, and apply Cauchy's formula for the oriented contour 
consisting of the segments $[-T,T], [T, T+ iy], [T+iy, -T +i y], [-T+iy, -T] $
\begin{eqnarray}
\int_{-T}^T e^{-itx} f(t)\,dt + 
\int_0^y e^{-i(T + ih) x} f(T+ih)\,dh
 & & \nonumber \\
 & & \hskip-70mm = \
\int_{-T}^T e^{-i(t + iy) x} f(t+iy)\,dt + 
\int_0^y e^{-i(-T + ih) x} f(-T+ih)\,dh.
\end{eqnarray}
Here, the two integrals taken over the interval $[0,y]$ are vanishing as 
$T \rightarrow \infty$. To prove this, first let us note that the functions 
$$
q_{k,h}(x) = e^{-hx} q_k(x), \quad x \in \R \ \ (k = 1,2),
$$
are integrable for every $h \in \R$ and have the Fourier transform 
$$
\widehat q_{k,h}(t) = \int_{-\infty}^\infty e^{itx}\,q_{k,h}(x)\,dx =
\E\,e^{i(t + ih)\eta_k} = f_k(t + ih).
$$
We may therefore conclude by applying the Riemann-Lebesgue lemma that 
$f_k(t + ih) \rightarrow 0$ as $|t| \rightarrow \infty$. Moreover, 
this convergence is uniform over all $0 \leq h \leq y$, which is due 
to the assumption $b)$. Indeed, since the mapping $h \rightarrow q_{k,h}$ 
from $[0,h]$ to $L^1(\R)$ is continuous, for any $\ep > 0$, one can choose
the points $0 = h_0 < h_1 < \dots < h_N = y$ such that
$\|q_{k,h} - q_{k,h_j}\|_{L^1} < \ep$ for all $h \in [h_j,h_{j+1}]$,
$0 \leq j \leq N-1$. In particular, 
$\sup_t |\widehat q_{k,h}(t) - \widehat q_{k,h_j}(t)| < \ep$.
By the Riemann-Lebesgue lemma, for every $j$, there is $t_j > 0$ such that
$\sup_{|t| \geq t_j} |\widehat q_{k,h_j}(t)| < \ep$. As a consequence,
$$
\sup_{h \in [0,y]}\, \sup_{|t| \geq T} |f_k(t+ih)| < 2\ep, 
$$
by choosing $T = \max\{t_1,\dots,t_N\}$. A similar property holds true
for $g_k$, $k = 1,2$, and therefore for the characteristic function $f$ 
in (11.6), we get
$$
\sup_{h \in [0,y]}\, \sup_{|t| \geq T}  |f(t+ih)| \rightarrow 0 \quad
{\rm as} \ \ T \rightarrow \infty.
$$

As a result, in the limit as $T \rightarrow \infty$ the identity (11.8)
leads to the equivalent variant of (11.7),
$$
p(x) = \frac{e^{yx}}{2\pi}
\int_{-\infty}^\infty e^{-itx} f(t+iy)\,dt,
$$
which yields
\be
p(x) \leq \frac{e^{yx}}{2\pi}
\int_{-\infty}^\infty |f(t+iy)|\,dt.
\en

In the next step we need to estimate the integrand in (11.9).
In view of the bound
$$
|g_0(t+iy)| = |\E\,e^{ic_0\, (t+iy) \eta_0}| \leq 
\E\,e^{-c_0 y\, \eta_0} = g_0(iy),
$$
(11.6) gives
$$
|f(t+iy)| \leq g_0(iy)\, |g_1(t+iy)|\, |g_2(t+iy)|.
$$
Applying this in (11.9) and using Cauchy's inequality, we get
\bee
p(x)
 & \leq &
e^{yx} g_0(iy)\ \frac{1}{2\pi} \int_{-\infty}^\infty 
|g_1(t+iy)|\, |g_2(t+iy)|\,dt \\
 & & \hskip-15mm \leq \
 e^{yx} g_0(iy)\ \bigg(\frac{1}{2\pi} \int_{-\infty}^\infty 
|g_1(t+iy)|^2\,dt\bigg)^{1/2}
\bigg(\frac{1}{2\pi} \int_{-\infty}^\infty 
|g_2(t+iy)|^2\,dt\bigg)^{1/2} \\
 & & \hskip-15mm = \
\frac{e^{yx}  f_0(ic_0 y)}{2\pi \sqrt{c_1 c_2}}
\bigg(\int_{-\infty}^\infty 
|f_1(t+ic_1 y)|^2\,dt\bigg)^{1/2}
\bigg(\int_{-\infty}^\infty 
|f_2(t+ic_2 y)|^2\,dt\bigg)^{1/2}.
\ene
Applying the Plancherel theorem and using the pointwise
subgaussian bound in $b)$, we get
\bee
\frac{1}{2\pi} \int_{\infty}^\infty |f_k(t+ic_k y)|^2\,dt
 & = &
\int_{\infty}^\infty e^{-2c_k yx} q_k^2(x)\,dx \\
 & \leq &
M_k^2 \int_{\infty}^\infty e^{-2c_k yx} \, \varphi^2(x)\,dx
 \, = \, 
\frac{M_k^2}{2\sqrt{\pi}}\,e^{c_k^2 y^2}.
\ene
In addition, by the assumption $a)$,
$f_0(ic_0 y) = \E\,e^{-c_0 y\, \eta_0} \leq e^{c_0^2 y^2/2}$.
Combining these estimates, we arrive at
$$
p(x) \leq \frac{e^{yx}}{\sqrt{2c_1 c_2}}\,\frac{M_1 M_2}{\sqrt{2\pi}}\,
e^{(c_0^2 + c_1^2 + c_2^2)\, y^2/2}.
$$
It remains to choose $y = -x$ and recall the assumption
$c_0^2 + c_1^2 + c_2^2 = 1$.
\qed

We conclude this section with 
\vskip2mm
{\bf Proof of Corollary 11.2.} To apply Theorem 2.2, 
we only need to check that $X$ has a density $p(x)$ satisfying (11.3).
Let $q(x)$ denote the common density of $\xi_k$, which is 
supposed to be bounded and compactly supported.
Without loss of generality, let $a_1 \geq a_2 \geq \dots \geq 0$.

{\it Case} 1: $a_1 = 1$ and $a_n = 0$ for all $n \geq 2$.
Then $p=q$, so that $p(x) \leq M_1 \varphi(x)$ a.e. 
for some constant $M_1 \geq 1$.

{\it Case} 2: $a_2 > 0$. Then 
$X = c_0 \eta_0 + c_1 \eta_1 + c_2 \eta_2$, where
$$
c_0 \eta_0 = \sum_{n = 3}^\infty a_n \xi_n, \quad
\eta_1 = \xi_1, \ \eta_2 = \xi_2, \ c_1 = a_1, \ c_2 = a_2.
$$
If $a_3 > 0$, then $c_0 = \sqrt{1 - a_1^2 - a_2^2}$,
so, $\eta_0$ is well-defined, strictly-subgaussian, and
has variance one. Otherwise, we may put $c_0 \eta_0 = 0$.
By Lemma 11.3, the relation $p(x) \leq M \varphi(x)$ a.e.
holds true with constant $M = \frac{1}{\sqrt{2a_1 a_2}}\,M_1^2$,
thus proving (11.3).
\qed

\vskip4mm
{\bf Acknowledgement.}
The research has been supported by the NSF grant DMS-2154001 and 
the GRF -- SFB 1283/2 2021 -- 317210226.


\vskip5mm
\begin{thebibliography}{BH3}
\itemsep=1pt

\bibitem{A-M-N}
Arbel, J.; Marchal, O.; Nguyen, H. D. On strict sub-Gaussianity, optimal 
proxy variance and symmetry for bounded random variables. 
ESAIM Probab. Stat. 24 (2020), 39--55.

\bibitem{A-B-B-N1} 
Artstein, S.; Ball, K. M.; Barthe, F.; Naor, A. Solution of Shannon's problem 
on the monotonicity of entropy. J. Amer. Math. Soc. 17 (2004), no. 4, 975--982. 

\bibitem{A-B-B-N2} 
Artstein, S.; Ball, K. M.; Barthe, F.; Naor, A. On the rate of convergence in 
the entropic central limit theorem. Probab. Theory Related Fields 129 (2004), 
no. 3, 381--390. 

\bibitem{Bar} 
Barron, A. R. Entropy and the central limit theorem. 
Ann. Probab. 14 (1986), no. 1, 336--342. 

\bibitem{B}
Bobkov, S. G. Upper bounds for Fisher information. 
Electron. J. Probab. 27 (2022), Paper No. 115, 44 pp.

\bibitem{B-C}
Bobkov, S. G.; Chistyakov, G. P. Bounds for the maximum of the density of the 
sum of independent random variables. (Russian) Zap. Nauchn. Sem. S.-Peterburg. 
Otdel. Mat. Inst. Steklov. (POMI) 408, Veroyatnost i Statistika. 18 (2012), 62--73, 324; 
transl. in J. Math. Sci. (N.Y.) 199 (2014), no. 2, 100--106.

\bibitem{B-C-G1}
Bobkov, S. G.; Chistyakov, G. P.; G\"otze, F.  Rate of convergence and 
Edgeworth-type expansion in the entropic central limit theorem. 
Ann. Probab. 41 (2013), no. 4, 2479--2512.

\bibitem{B-C-G2}
Bobkov, S. G.; Chistyakov, G. P.; G\"otze, F. Berry-Esseen bounds in 
the entropic central limit theorem. Probab. Theory Related Fields 
159 (2014), no. 3-4, 435--478.

\bibitem{B-C-G3}
Bobkov, S. G.; Chistyakov, G. P.; G\"otze, F. R\'enyi divergence and the central limit 
theorem. Ann. Probab. 47 (2019), no. 1, 270--323. 

\bibitem{B-C-G4}
Bobkov, S. G.; Chistyakov, G. P.; G\"otze, F. Richter's local limit theorem, its 
refinement, and related results. Lithuanian J. Math. 63 (2023), no. 2, 138--160.

\bibitem{B-C-G5}
Bobkov, S. G.; Chistyakov, G. P.; G\"otze, F. Strictly subgaussian distributions.
Preprint (2023). To appear in: Electron. J. Probab.

\bibitem{B-G}
Bobkov, S. G.; G\"otze, F. Exponential integrability and transportation cost 
related to logarithmic Sobolev inequalities. J. Funct. Anal. 163 (1999), 
no. 1, 1--28. 

\bibitem{Bu-K1}
Buldygin, V. V.; Kozachenko, Yu. V.  Sub-Gaussian random variables. (Russian) 
Ukrain. Mat. Zh. 32 (1980), no. 6, 723--730.

\bibitem{Bu-K2}
Buldygin, V. V.; Kozachenko, Yu. V. Metric characterization of random variables 
and random processes. Translated from the 1998 Russian original by V. Zaiats. 
Transl. Math. Monogr., 188
American Mathematical Society, Providence, RI, 2000. xii+257 pp.

\bibitem{D} 
Daniels, H. E. Saddlepoint approximations in statistics. 
Ann. Math. Statist. 25 (1954), 631--650.

\bibitem{D-C-T} 
Dembo, A., Cover, T. M., Thomas, J. A. Information-theoretic inequalities. 
IEEE Trans. Inform. Theory, 37 (1991), no. 6, 1501--1518.

\bibitem{E}
Esscher, F. On the probability function in the collective theory of risk. 
Skandinavisk Aktuarietidskrift. 15 (3) (1932), 175--195.

\bibitem{G-H}
Guionnet, A.; Husson, J. Large deviations for the largest eigenvalue of 
Rademacher matrices. Ann. Probab. 48 (2020), no. 3, 1436--1465.

\bibitem{H-N-T}
Havrilla, A.; Nayar, P.; Tkocz, T. Khinchin-type inequalities via Hadamard's 
factorisation. Int. Math. Res. Not., no. 3 (2023), 2429--2445.

\bibitem{I-L}
Ibragimov, I. A.; Linnik, Yu. V. Independent and stationary sequences of 
random variables. With a supplementary chapter by I. A. Ibragimov and V. V. Petrov.
Translation from the Russian edited by J. F. C. Kingman. 
Wolters-Noordhoff Publishing, Groningen, 1971. 443 pp. 

\bibitem{J}
Johnson, O. Information theory and the central limit theorem. Imperial College Press, 
London, 2004. xiv+209 pp.

\bibitem{K}
Khinchin, A. I. Mathematical Foundations of Statistical Mechanics.
Translated by G. Gamow. Dover Publications, Inc., New York, N.Y., 1949. viii+179 pp.

\bibitem{LC}
Le Cam, L. Asymptotic methods in statistical decision theory. 
Springer Series in Statistics. Springer-Verlag, New York, 1986. xxvi+742 pp.

\bibitem{L-Y}
Lee, T.-D.; Yang, C.-N. Statistical theory of equations of state and phase transitions. II. 
Lattice gas and Ising model. Phys. Rev. 87 (3), 410 (1952).

\bibitem{M-B}
Madiman, M.; Barron, A. Generalized entropy power inequalities and 
monotonicity properties of information. 
IEEE Trans. Inform. Theory 53 (2007), no. 7, 2317--2329. 

\bibitem{N1}
Newman, C. M. Inequalities for Ising models and field theories which obey the Lee-Yang 
theorem. Comm. Math. Phys. 41 (1975), 1--9.

\bibitem{N2}
Newman, C. M. Moment inequalities for ferromagnetic Gibbs distributions.
J. Mathematical Phys. 16 (1975), no. 9, 1956--1959.

\bibitem{N3}
Newman, C. M. An extension of Khintchine's inequality. Bull. 
Amer. Math. Soc. 81 (1975), no. 5, 913--915.

\bibitem{N-W}
Newman, C.; Wu, W. Lee-Yang property and Gaussian multiplicative chaos.
Comm. Math. Phys. 369 (2019), no. 1 , 153--170.

\bibitem{Pe1}
Petrov, V. V. Local limit theorems for sums of independent random variables. 
            (Russian) Teor. Verojatnost. i Primenen. 9 (1964), 343--352. 

\bibitem{Pe2}
Petrov, V. V. Sums of independent random variables. Translated from 
the Russian by A. A. Brown. Ergebnisse der Mathematik und ihrer 
Grenzgebiete, Band 82. Springer-Verlag, New York-Heidelberg, 1975. x+346 pp.

\bibitem{Pr}
Prokhorov, Yu. V. A local theorem for densities. (Russian) 
Doklady Akad. Nauk SSSR (N.S.) 83 (1952), 797--800. 

\bibitem{Ri}
Richter, W. Lokale Grenzwerts\"atze f\"ur grosse Abweichungen. (Russian) 
Teor. Veroyatnost i Primenen. 2 (1957), 214--229. 

\bibitem{S}
Shiryaev, A. N. Probability. Translated from the first (1980) Russian edition 
by R. P. Boas. Second edition. Graduate Texts in Mathematics, 95. 
Springer-Verlag, New York, 1996. xvi+623 pp.

\bibitem{VE-H}
van Erven, T.; Harremo\"es, P. R\'enyi divergence and Kullback-Leibler divergence.
IEEE Trans. Inform. Theory 60 (2014), no. 7, 3797--3820.

\end{thebibliography}
\end{document}